\documentclass[11pt,
]{article}

\usepackage{amssymb,amstext,amsthm,latexsym}
\usepackage{amsmath}
\usepackage{amsfonts}
\usepackage{mathrsfs}
\usepackage{mparhack}
\usepackage[bbgreekl]{mathbbol}

\newcommand{\nek}{\newcommand}
\nek{\renek}{\renewcommand}

\makeatletter
\if@twoside%
\else\fi%
\makeatother

\nek{\parf}{\subsection}
\nek{\punk}{\subsection}

\nek{\skr}{\mathscr}
\nek{\cP}{\skr P}

\nek{\vyk} [1] {}

\nek{\imar}[1] 
{}
\vyk{
{\marginpar[
\flushright\footnotesize%
$\mtho\longrightarrow$\\%
\vspace{-1ex}{#1}\vspace*{1ex}]%
{
\flushleft\footnotesize%
$\mtho\longleftarrow$\\%
\vspace{-1ex}{#1}\vspace*{1ex}}%
}%
}

\nek{\imae} {\imar}

\nek{\itsep}{\itemsep=0.3ex plus 0.1ex minus 0.1ex}

\nek{\tenu}[1]{
\def\theenumi{#1}
\itsep
}

\nek{\tenud}[1]{
\def\theenumi{#1}
\itsep
}

\nek{\cenu}{\tenud{{%
$\mtho\arabic{enumi}^\circ$%
}}}

\theoremstyle{plain}

\newtheorem{theore}             {Theorem} 
\newtheorem{corollar}  [theore]{Corollary}
\newtheorem{propo}     [theore]{Proposition}
\newtheorem{lemm}      [theore]{Lemma}
\newtheorem{lemt}      {Lemma}[theore]
\newtheorem{cla}       [theore]{Claim} 
\newtheorem{clt}       {Claim} [theore]

\theoremstyle{definition}
\newtheorem{defn}      [theore]{Definition}
\newtheorem{rem} [theore]   {Remark} 
\newtheorem*{solu}{Solution}      
\newtheorem*{probl}{Problem}      
\newtheorem{exz} {Example} 
\newtheorem*{prF}{Proof}               

\nek{\thsp}{\hspace{0.1ex plus \mathsurround}}

\nek{\bex}{\begin{exz}}
\nek{\eex}{\end{exz}}
\nek{\bprl}{\begin{probl}}
\nek{\eprl}{\end{probl}}
\nek{\bsol}{\begin{solu}}
\nek{\esol}{\end{solu}}
\nek{\bpro}{\begin{propo}}
\nek{\epro}{\end{propo}}
\nek{\bcor}{\begin{corollar}}
\nek{\ecor}{\end{corollar}}
\nek{\bdf} {\begin{defn}} 
\nek{\eDf} {\end{defn}}
\nek{\edf} {\qed\end{defn}}
\nek{\edF} {\end{defn}}
\nek{\ble} {\begin{lemm}}
\nek{\ele} {\end{lemm}}
\nek{\blt} {\begin{lemt}}
\nek{\elt} {\end{lemt}}
\nek{\bre}{\begin{rem}}
\nek{\ere}{\end{rem}}
\nek{\bte} {\begin{theore}}
\nek{\ete} {\end{theore}}
\nek{\bpf} {\begin{prF}} 
\nek{\epf} {\qed\end{prF}} 
\nek{\ePf} {\end{prF}} 
\nek{\qeD} {\qed}
\nek{\qeDD} [1] 
{\hfill\hbox{\qed~({\small #1\/}\hspace{0.1ex})}}
\nek{\epF} [1] {\qeDD{#1}\end{prF}} 

\nek{\ben}{\begin{enumerate}\itsep}
\nek{\een}{\end{enumerate}}
\nek{\bit}{\begin{itemize}\itsep}
\nek{\eit}{\end{itemize}}
\nek{\bde}{\begin{description}\itsep}
\nek{\ede}{\end{description}}
\nek{\bay}{\begin{array}}
\nek{\eay}{\end{array}}
\nek{\bce}{\begin{center}}
\nek{\ece}{\end{center}}
\nek{\bqo}{\begin{quotation}\noi}
\nek{\eqo}{\end{quotation}}
\nek{\ZFC}{{\bf ZFC}}
\nek{\ZF}{{\bf ZF}}

\nek{\iesp}{\hspace{0.3ex}}
\nek{\resp}{\hspace{0.25ex}}
\nek{\ie} {\hbox{\sl i.\iesp e.}}
\nek{\eg} {\hbox{\sl e.\iesp g.}}
\nek{\noo}
{\hbox{\sl w.\iesp l.\iesp o.\iesp g.}}
\nek{\vrt} {\hbox{w.\iesp r.\iesp t.}}
\nek{\ddd}[1]{$\mtho\hspace{0.2ex}{#1}$-\hspace{0.0ex}}
\nek{\dd}{\ddd}
\nek{\ran}  {\mathop{\tt ran}}
\nek{\otp} [1] {\mathop{\rm otp}(#1)}
\nek{\dom}  {\mathop{\tt dom}}
\nek{\cof}  {\mathop{\tt cof}}
\nek{\tsup} {\mathop{\tt sup}}
\nek{\tmax} {\mathop{\tt max}}
\nek{\tmin} {\mathop{\tt min}}
\nek{\tlim} {\mathop{\tt lim\hspace{0.3ex}}}

\nek{\al} {\alpha}
\nek{\ba} {\beta}
\nek{\ga} {\gamma}
\nek{\da} {\delta}
\nek{\za} {\zeta}
\nek{\ka} {\kappa}
\nek{\La}{\Lambda}
\nek{\la}{\lambda}
\nek{\sg} {\sigma}
\nek{\Sg} {\Sigma}
\nek{\ve}{\varepsilon}
\nek{\vpi}{\varphi}
\nek{\vt} {\vartheta}
\nek{\om} {\omega}
\nek{\Om} {\Omega}
\nek{\lom}{^{<\om}}
\nek{\omi} {\om_1}
\nek{\ali} {\aleph_1}
\nek{\alo} {{\aleph_0}}

\nek{\fs}[2]{{\hspace*{0.3ex}\boldsymbol\Sigma}^{#1}_{#2}}
\nek{\fp}[2]{{\boldsymbol\Pi}^{#1}_{#2}}
\nek{\fd}[2]{{\boldsymbol\Delta}^{#1}_{#2}}

\nek{\id}[2]{{\varDelta}^{#1}_{#2}}
\nek{\ip}[2]{{\varPi}^{#1}_{#2}}
\nek{\is}[2]{{\varSigma}^{#1}_{#2}}

\nek{\BBB}{\hspace{0.01ex}}
\nek{\dR}{{\BBB{\mathbb R}\BBB}}
\nek{\dP}{{\BBB{\mathbb P}\BBB}}

\nek{\bn}{\om^\om}
\nek{\sus} {\mathopen{\exists\hspace{0.35ex}}}
\nek{\kaz} {\mathopen{\forall\hspace{0.35ex}}}
\nek{\imp} {\Longrightarrow} 
\nek{\mpi} {\Longleftarrow} 
\nek{\eqv} {\Longleftrightarrow} 
\nek{\leqv} {\;\eqv\;} 
\nek{\limp} {\;\imp\;} 
\nek{\ti}  {\times} 
\nek{\sq}  {\subseteq}
\nek{\su}  {\subset}
\nek{\sneq}{\subsetneqq}
\nek{\we}  {{\mathbin{\hspace{0.15ex}^\wedge}}}
\nek{\obr} {^{-1}}
\nek{\dif} {\smallsetminus}
\nek{\res} {\mathbin{\restriction}}
\nek{\lef} {\preccurlyeq}
\nek{\gef} {\succcurlyeq}
\nek{\pu}  {\varnothing}
\nek{\iy}  {\infty}
\nek{\piy} {+\iy}
\nek{\nin} {\not\in}
\nek{\onto}{\stackrel{\text{\rm onto}}{\longrightarrow}}
\nek{\ang} [1] {\langle #1\rangle}
\nek{\stk} [2] {\ang{#1\hspace{0.3ex};\hspace{0.1ex}#2}}
\nek{\sis} [2] {\ans{#1}_{#2}}
\nek{\ans} [1] {\{\hspace{0.01ex}#1\hspace{0.01ex}\}}

\nek{\zz} {\linebreak[0]} 
\nek{\ens} [2] {\ans{{#1\hspace{0.5ex}{:}}\zz\hspace{0.5ex}#2}}
\nek{\itla} {\item\label}

\nek{\ubf}{\fontseries{b}\selectfont}
\nek{\TS}{\textstyle}
\nek{\yo} {,\linebreak[0]}
\nek{\yi} {\hspace{0.2ex},\linebreak[0]\hspace{0.2ex}}
\nek{\yd} {\hspace{0.2ex},\linebreak[0]\:}
\nek{\yt} {\hspace*{0.2ex},\linebreak[0]\;} 
\nek{\lex} {<_{\text{\tt lex}}}
\nek{\lexe} {\leqslant_{\text{\tt lex}}}

\nek{\snos} [1] {\,\footnote{\hspace*{2pt}#1}}
\nek{\renu}{\tenu{{\rm(\roman{enumi})}}}
\nek{\fenu}{\tenu{{\rm(\fnsymbol{enumi})}}}
\nek{\Renu}{\tenu{{\rm(\Roman{enumi})}}}
\nek{\Renup}{\tenu{{\rm(\Roman{enumi}$'$)}}}
\nek{\rit} [1] {{\it#1\/}}
\nek{\lam} [1] {\label{#1}\hspace*{-3pt}\imar{#1}}%
\nek{\las} [1] {\label{#1}\imar{#1}}%

\nek{\atc}{\addtocounter{enumi}{1}}
\nek{\atcm}{\addtocounter{enumi}{-1}}

\nek{\wo} {\text{\ubf WO}}
\nek{\bez} {\dif}
\nek{\AC}   {{\text{\bf AC}}}
\nek{\DC}   {{\text{\bf DC}}}
\nek{\od} {\text{OD}}
\nek{\rod} {\text{ROD}}
\nek{\ROD}  {{\text{\bf{ROD}}}}
\nek{\hrod}{{\text{\bf HROD}}}
\nek{\OD} {\mathbf P}

\nek{\rH} {\text{\ubf H}}
\nek{\rL} {\text{\ubf L}}
\nek{\rV} {\text{\ubf V}}
\nek{\rU} {\text{\ubf U}}

\nek{\led} {\leq^\ast}
\nek{\ld} {<^\ast}

\nek{\pqo} {PQO}

\nek{\Ord}  {\mathop{\tt Ord}}
\nek{\supp}  {\mathop{\tt supp}}
\nek{\card}  {\mathop{\tt card}}

\nek{\np}{\newpage}

\nek{\mtho}{\mathsurround=0mm}
\nek{\msur}{\hspace{-1\mathsurround}}
\nek{\noi}{\noindent}
\nek{\vom}{\vspace{1mm}}
\nek{\vim}{\vspace{-1mm}}
\nek{\vtm}{\vspace{2mm}}

\nek{\eqr} {equivalence relation}

\nek{\rE} {\mathrel{\mathsf E}}

\nek{\ek}[2] {[#1]_{{#2}}}

\nek{\eke}[1] {\ek{#1}{\rE}}

\nek{\qand}{\quad\text{and}\quad}

\nek{\dX}{{\BBB{\mathbb X}\BBB}}
\nek{\dN}{{\BBB{\mathbb N}\BBB}}
\nek{\dB}{{\BBB{\mathbb B}\BBB}}
\nek{\dF}{{\BBB{\mathbb F}\BBB}}
\nek{\fun}{{\BBB{\mathbb {Fun}}\BBB}}

\nek{\doP}  [1] {{#1}^\complement}

\nek{\curle}{\preccurlyeq}
\nek{\cle}{\curle}
\nek{\cl} {\prec}
\nek{\ncle}{\not\cle}
\nek{\ncl}{\not\cl}

\nek{\gh}{\mathbb P}
\nek{\ghd}{\gh^2}

\nek{\fdt}{\hbox{\raisebox{-0.25ex}{\LARGE\bf.}}}
\nek{\bdot}[1] 
{\raisebox{-0.07ex}{\mtho$\stackrel{\fdt}{#1}$}}

\nek{\ile} {_{\text{\tt le}}}
\nek{\iri} {_{\text{\tt ri}}}

\nek{\dox}{\bdot{\boldsymbol x}}
\nek{\doy}{\bdot{\boldsymbol y}}

\nek{\doxl}{\dox\ile}
\nek{\doxr}{\dox\iri}
\nek{\doyl}{\doy\ile}
\nek{\doyr}{\doy\iri}

\nek{\dotx}{\check x}
\nek{\dotY}{\check Y}

\nek{\dn}{2^\om}
\nek{\dpe} {\mathord{\drof\gh\rE}}
\nek{\dpew}{\mathord{\drow{\gh}\rE}}

\nek{\trof} [3] {{#1}\ti_{#2}{#3}}

\nek{\drof} [2] {{#1}\ti_{#2}{#1}}
\nek{\drow} [2] {{#1}\ti_{#2}^{\text{weak}}{#1}}

\nek{\ck} {\om_1^{\text{\sc ck}}}

\nek{\Eo}  {\rE_{\text{\sf0}}}
\nek{\Fo}  {\rF_{\text{\sf0}}}
\nek{\nE}  {\mathbin{{\not\hspace{-0.35ex}\sf E}}}

\nek{\bcl} {\begin{cla}}
\nek{\ecl} {\end{cla}}
\nek{\bct} {\begin{clt}}
\nek{\ect} {\end{clt}}

\nek{\bV}{{\mathbf V}}
\nek{\bL}{{\bf L}}

\nek{\gM} {\mathfrak M}
\nek{\gN} {\mathfrak N}
\nek{\cM} {\skr M}
\nek{\cN} {{\skr N}}
\nek{\mN} {\boldsymbol N}

\nek{\ccs} {}
\nek{\cF}{{\ccs{\skr F}\ccs}}
\nek{\cS}{{\ccs{\skr S}\ccs}}
\nek{\cD}{{\ccs{\skr D}\ccs}}

\nek{\cf} [1] {\cF_{#1}}
\nek{\cfx} [1] {\cF_#1}
\nek{\cfd} [2] {\cF_{#1}(#2)}

\nek{\dfx} [2] {\cF_#1^#2}

\nek{\etc} {{\sl etc}}

\nek{\bus}{\begin{equation}}   
\nek{\eus}{\end{equation}}

\nek{\pp} [2] {\gh^{#1}_{#2}}

\nek{\dpd} [1] {\gh\ti_{#1}\gh}
\nek{\dpx} [1] {\gh\ti_{\rE_{#1}}\gh}
\nek{\dpp} {\dpd{}}

\nek{\dpw} {\gh(W)}
\nek{\dpwe} {\gh(W)\ti_{\rE}\gh(W)}

\nek{\odw} {\odsw{W}}
\nek{\odwp} {\odsw{W'}}
\nek{\odwe} {\odw\ti_{\rE}\odw}
\nek{\odwep} {\odwp\ti_{\rE}\odwp}

\nek{\ups} {\varUpsilon}

\nek{\apr} {\approx}
\nek{\napr}{\not\apr}

\nek{\gp}{\mathfrak p}
\renek{\gp}{\mathbb p}

\nek{\aenu}{\tenu{{\rm(\arabic{enumi})}}}
\nek{\aenup}{\tenu{{\rm(\arabic{enumi}$'$)}}}
\nek{\Aenu}{\tenu{{\rm(\Alph{enumi})}}}
 
\nek{\doxlp}{\dox{}'_{\tt le}}
\nek{\doxrp}{\dox{}'_{\tt ri}}

\nek{\bfit}{\bfseries\itshape}

\nek{\esn}{\rS_{\ans{1/n}}}
\nek{\rS}  {\mathbin{\sf S}}

\nek{\wh}{\widehat}
\nek{\wY}{\widehat Y}
\renek{\wY}{C}
\nek{\wX}{C}

\nek{\lr}{LR}
\nek{\rl}{RL}

\nek{\vx} {\vec x}

\nek{\Pa}{P^*}
\nek{\Xa}{X^*}
\nek{\Ua}{U^*}
\nek{\Va}{V^*}
\nek{\Do}{D_0}

\nek{\rEF} {\mathrel{\mathsf E_\cF}}
\nek{\refx} [1] {\mathrel{\mathsf E_{\cfx#1}}}

\nek{\osm} {\dd\Omega SM}
\nek{\odk} {\dd\od1st-countable\/}
\nek{\psur}{\hspace{1\mathsurround}}
\nek{\pwod}  [1] {\cP_{\text{\tt OD}}(#1)}
\nek{\pws}  [1] {\cP(#1)}
\nek{\odi} {\OD^*}
\nek{\koh} {\text{\sc Coh}}
\nek{\nse} {\om\lom}
\nek{\bse} {2\lom}
\nek{\ibn} [1] {\skr N_{#1}}

\nek{\ilom}{\omi\lom}

\nek{\vT} {\Theta}

\nek{\col} [1] {\text{\ubf Coll}(\om,#1)}
\nek{\coll} [1] {\col{{<\hspace*{0.1ex}}#1}}

\nek{\msl} {\mathbin{<_{\text{lex}}}} 
\nek{\Lom} {^{<\Om}} 

\nek{\fh} {{\mathfrak h}}
\nek{\ghx} {\text{\ubf h}} 
    \nek{\ghi} {\ghx^{\omi}}
    \nek{\gho} {\fh}
    \nek{\ghp} {\fh'}

\nek{\dodg} {\od-generic}
\nek{\pge} {\dd\OD generic}
\nek{\ege} {\dd{(\spe)}generic}
\nek{\ega} [1] {\dd{(\spa{#1})}generic}

\nek{\spe} {\mathord{\drof\OD\rE}}
\nek{\spa} [1] {\mathord{\drof\OD{\rE_{#1}}}}

\nek{\spei} 
{\mathord{\drof{\odi\hspace*{-0.3ex}}\rE\hspace*{0.3ex}}}
\nek{\spai} [1] 
{\mathord{\drof{\odi\hspace*{-0.3ex}}{\rE_{#1}}\hspace*{0.3ex}}}

\nek{\xle} {x\ile}
\nek{\xri} {x\iri}

\nek{\odf} {\od-force}
\nek{\pfo} {\dd\OD force}
\nek{\efo} {\dd{(\spe)}force}
\nek{\efa} [1] {\dd{(\spa{#1})}force}

\nek{\odsw} [1] {\OD_{\sq#1}}
\nek{\odd}   {\OD^{(2)}}

\nek{\pe} [1] {{\text{\ubf Perf}\hspace*{0.1ex}}_{#1}}
\nek{\cs} [2] {(#1)_{#2}}

\nek{\gpe} [1] {{\gh}_{#1}}
\nek{\gped} {\gpe2\ti\gpe2}
\nek{\gpei} {\gpe2\tj 1\gpe2}

\nek{\doxi}{\dox_1}
\nek{\doxd}{\dox_2}

\nek{\ilei} {_{\text{\tt le},1}}
\nek{\irii} {_{\text{\tt ri},1}}
\nek{\iled} {_{\text{\tt le},2}}
\nek{\irid} {_{\text{\tt ri},2}}

\nek{\doxli}{\dox\ilei}
\nek{\doxri}{\dox\irii}
\nek{\doxld}{\dox\iled}
\nek{\doxrd}{\dox\irid}
                            
\nek{\tj} [1] 
{\mathop{{\hspace*{0.2ex}\ti\hspace*{-0.3ex}}_%
{\text{$\scriptstyle#1$}}}}

\nek{\comp} [1] {\text{\tt Comp}_{#1}}

\nek{\rdi} {\la}
\nek{\mekt} {\mek_T}
\nek{\meks} {\mek_S}
\nek{\mek}  {\preccurlyeq}
\nek{\bbg}{\hspace{0.1ex}} 
\nek{\pr} {{\bbg{\fraK p}\bbg}\hspace{1pt}}
\nek{\fraK}{\mathfrak}

\nek{\dop} [1] {\mathopen\complement{\hspace{0.1ex}{#1}}}
\nek{\klst}  {{\kL[S,T]}}

\nek{\rdii}{(\rdi\pone)}
\nek{\pone}{\hspace{-0.4ex}+\hspace{-0.4ex}1}
\nek{\hop} {\hbox{h.\dosp o.\dosp p.}}
\nek{\dosp}{\hspace{0.4ex}}

\nek{\kL} {\rL}
\nek{\kV} {\rV}

\nek{\mel} {\mathbin{\leq_{\text{lex}}}}
\nek{\ees} {\approx_S}
\nek{\meo} {\mathbin{\leq_0}}
\nek{\nEo} {\mathbin{\not{\Eo}}}
\nek{\nmeks}{\not\meks}

\nek{\ult} [2] {\text{\ubf Ult}_{#1}(#2)}
\nek{\ugn} {\ult G\gN}
\nek{\zfc} {\ZFC}

\nek{\hkp} {\rH_{\ka'}}
\nek{\er} [1] {\mathrel{\rE_#1}}
\nek{\ear} [2] {\mathrel{\rE_#1^#2}}

\nek{\aprs} {\mathrel{\apr_S}}
\nek{\cles} {\mathrel{\cle_S}}
\nek{\cls} {\mathrel{\cl_S}}

\nek{\bda} {\boldsymbol\da}
\nek{\kap} {\ka^+}

\begin{document}

\title
{On countable cofinality and decomposition 
of definable thin orderings}

\author 
{Vladimir~Kanovei\thanks{IITP RAS and MIIT,
Moscow, Russia, \ {\tt kanovei@googlemail.com}. 
Partial support of RFFI grant 13-01-00006 acknowledged. 
--- \it Contact author.}  
\and
Vassily~Lyubetsky\thanks{IITP RAS,
Moscow, Russia, \ {\tt lyubetsk@iitp.ru}} 
}

\date 
{\today}

\maketitle

\begin{abstract}
We prove that in some cases definable thin sets 
(including chains) of Borel
partial orderings are necessarily countably cofinal.
This includes the following cases:
analytic thin sets,
\ROD\ thin sets in the Solovay model,
and $\fs12$ thin sets in the assumption that
$\omi^{\rL[x]}<\omi$ for all reals $x$.
We also prove that definable thin wellorderings admit 
partitions into definable chains in the Solovay model.
\end{abstract}

\parf{Introduction}
\las{intro}
 
Studies of maximal chains in partially ordered sets go back 
to as early as Hausdorff \cite{h07,h09}, where this issue 
appeared in connection with Du Bois Reymond's investigations 
of orders of infinity.
Using the axiom of choice, Hausdorff proved the existence of
maximal chains (which he called \rit{pantachies}) 
in any partial ordering.
On the other hand, Hausdorff clearly understood the difference
between such a pure existence proof and an actual construction
of a maximal chain --- see \eg\ \cite[p.~110]{h07} or
comments in \cite{fi} --- which we would understand nowadays 
as the existence of \rit{definable} maximal chains.

The following theorem present three cases in which all 
linear, and even thin suborders 
of Borel PQOs are necessarily countably cofinal. 

\vyk{
\bte
\lam{mt}
If\/ $\cle$ is a Borel PQO on a (Borel) set\/ 
$D=\dom{(\cle)}\sq\bn$,  
$X\sq D$, and\/ ${\cle}\res X$ is a thin quasi-ordering\/
then%
$:$ 
\ben
\renu
\itla{mt1}
if\/ $X$ is a\/ {$\fs11$} set then 
there is no\/ \dd\cle chains in\/ $X$ of uncountable 
cofinality, and\/ $X$ itself is countably\/ \dd\cle cofinal,   

\itla{mt2}
if\/ $X$ is a\/ \ROD\ set in the Solovay model, then
there is no cofinal\/ \dd\cle chains in\/ $X$ of uncountable 
cofinality,

\itla{mt3}
if\/ $X$ is a\/ $\fs12$ set, and\/
$\omi^{\rL[r]}<\omi$ for every real\/ $r$, then
there is no cofinal\/ \dd\cle chains in\/ $X$ of uncountable 
cofinality. 
\een
Therefore, if, in addition, it is known that\/ $\stk D\le$
does not have cofinal countable chains, then
in all three cases\/ $X$ is not cofinal in\/ $D$.
\ete
} 

\bte
\lam{mt}
If\/ $\cle$ is a Borel PQO on a (Borel) set\/ 
$D=\dom{(\cle)}\sq\bn$,  
$X\sq D$, and\/ ${\cle}\res X$ is a thin quasi-ordering\/
then\/ $\ang{X;{\cle}}$ is countably cofinal in each of the 
following three cases$:$ 
\ben
\renu
\itla{mt1}
if\/ $X$ is a\/ {$\fs11$} set, and in this case moreover   
there is no\/ \dd\cle chains in\/ $X$ of uncountable 
cofinality,   

\itla{mt2}
if\/ $X$ is a\/ \ROD\ set in the Solovay model, 

\itla{mt3}
if\/ $X$ is a\/ $\fs12$ set, and\/
$\omi^{\rL[r]}<\omi$ for every real\/ $r$. 
\een
Therefore, if, in addition, it is known that\/ $\stk D\le$
is not countably cofinal, then
in all three cases\/ $X$ is not cofinal in\/ $D$.
\ete

The additional condition in the theorem,
of the uncountable cofinality,
holds for many partial orders of interest, \eg,
the eventual domination order on sets like $\om^\om$ or
$\dR^\om,$ or the \rit{rate of growth order} defined
on $\dR^\om$ by
$
x<_{\text{\sc rg}}y
\quad\text{iff}\quad
\lim_{n\to\infty}\frac{y(n)}{x(n)}=\infty 
$
(see a review in \cite{hos}).
Needless to say that chains, gaps, and similar structures
related to these or similar orderings have been subject of
extended studies, of which we mention \cite{aq,fa1,to:gap2,hom}
among those in which the definability aspect is considered.

Part \ref{mt1} of the theorem is proved 
in Section~\ref{sek11} by reduction 
to a result (Theorem~\ref{hmst} below) which 
extends a 
theorem in \cite{hms} to the case of $\fs11$ 
suborders of a background Borel PQO as in \ref{mt1}.
Part \ref{mt2} is already known from \cite{kl} in the case 
of linear, rather than thin, $\ROD$ suborders, but we
present here (Section~\ref{sbSM}) 
an essentially simplified proof. 
Part \ref{mt3} is proved in Section~\ref{sbS12} by a reference 
to part \ref{mt2} and a sequence of absoluteness arguments. 

It is a challenging question to figure out 
whether claims \ref{mt2} and \ref{mt3} of Theorem \ref{mt} 
remain true in stronger forms similar to the ``moreover'' 
form of claim \ref{mt1}.
The answer is pretty simple in the affirmative provided we 
consider only accordingly definable 
(but not necessarily cofinal) \ddd\omi sequences in 
the given set $X$ --- that is to say, \ROD\ in claim \ref{mt2} 
and $\fs12$ in claim \ref{mt3}.

The next theorem (our second main result) 
extends a classical decomposition theorem in \cite{hms} 
to the case of definable sets in the Solovay model. 

\bte
[in the Solovay model]
\lam{t:sm}
Let\/ $\cle$ be an\/ \od\ PQO on\/ $\bn$, 
$\apr$ be the associated \eqr, 
and\/ $\Xa\sq \bn$ be an \od\ \dd\cle thin set.
Then\/ 
$\Xa$ is covered by the union of\/ 
all\/ $\od$ \dd\cle chains\/ $C\sq\bn\;.$ 

The same is true for any definability class\/ $\od(x)$, 
where\/ $x$ is a real. 
\ete

The proof of this theorem 
as given in Section~\ref{Td:sm} has a certain 
semblance of the proof of Theorem 5.1 in \cite{hms} in 
the context of its general combinatorial structure. 
Yet the proof includes some 
changes necessary since $\od$ sets in the Solovay model 
only partially resemble sets in $\id11$ and $\is11$.
In particular we'll have to establish some properties of 
the $\od$ forcing rather different from the properties 
of the Gandy -- Harrington forcing applied in \cite{hms}, 
and also prove a tricky compression lemma 
(Lemma~\ref{apal31}) 
in Sections \ref{smod} -- \ref{CL}.

\parf{Notation}
\las{not}

We proceed with notational remarks.

%
\bde
\item[\rm PQO, \it partial quasi-order\/$:$] \ 
reflexive ($x\le x$) and transitive in the domain;


\item[\rm LQO,  
\it linear quasi-order\,$:$] \  
PQO and $x\le y\lor y\le x$  
in the domain; 

\item[\rm LO, \it linear order\,$:$]  \ 
LQO and $x\le y\land y\le x\imp x=y$; 

\item[\it associated equivalence relation\,$:$]  \ 
$x\apr y$ iff $x\le y\land y\le x$.

\item[\it associated strict order\,$:$]  \ 
$x<y$ iff $x\le y\land y\not\le x$.
\ede
By default we consider only \rit{non-strict} orderings.
All cases of consideration of \rit{strict} PQOs
will be explicitly specified.
\bde
\item[\it strict \rm PQO\,$:$] \ irreflexive ($x\not<x$) and
transitive; 


\item[\it strict \rm LO\,$:$]  \ 
strict PQO and the trichotomy $\kaz x,y\:(x<y\lor y<x\lor x=y)$.

\item[\it \lr\ (left--right) order preserving map$:$]
any map $f:\stk{X}{\le}\to\stk{X'}{\le'}$ such that 
we have $x\le y\imp f(x)\le' f(y)$ for all $x,y\in\dom f$;

\item[\it \rl\ (right--left) order preserving map$:$]
a map $f:\stk{X}{\le}\to\stk{X'}{\le'}$ such that 
we have $x\le y\mpi f(x)\le' f(y)$ for all $x,y\in\dom f$;


\item[\it sub-order\,$:$] \ 
a restriction of the given PQO to a subset of
its domain.

\item[\it $\lex\yt\lexe\;:$]
the lexicographical LOs on sets of the form $2^\al\yd\al\in\Ord$, 
resp.\ strict and non-strict. 
\ede

Let $\stk P\le$ be a background PQO. 
A subset $Q\sq P$ is:

\bde
\item[\it cofinal in $P\,:$] iff \
$\kaz p\in P\:\sus q\in Q\:(p\le q)$;

\item[\it countably cofinal (in itself)$:$] iff there exists 
a countable set $Q'\sq Q$ cofinal in $Q$;

\item[\it a chain$:$]
iff it consists of 2wise \dd\le comparable elements, \ie, LQO;

\item[\it an antichain in $P\,:$]
iff it consists of 2wise \dd\le incomparable elements;

\item[\it a thin set$:$]
iff it contains no perfect \dd\le antichains.
\ede


Finally if $\rE$ is an \eqr\ then let
$$
\bay{rcll}
\eke x&=&\ens{y\in\dom\rE}{x\rE y} &
\text{(the \dd\rE\rit{class} of $x\in\dom\rE$),}\\[1ex] 
\eke X&=&\textstyle\bigcup_{x\in X}\eke x & 
\text{(the \dd\rE\rit{saturation} of $X\sq\dom\rE$).}
\eay
$$


\punk{Analytic thin subsets}
\las{sek11}

In this Section, we prove Theorem \ref{mt}\ref{mt1} by reference 
to the following background result: 

\bte
[proved in \cite{kanHMS}]
\lam{hmst}
Let\/ $\cle$ be a\/ $\id11$ PQO on\/ $\bn$, 
$\apr$ be the associated \eqr, 
and\/ $\Xa\sq \bn$ be a\/ $\is11$ set.
Then\/ 
\ben
\Renu
\itla{hmst1} 
if\/ $\Xa$ is\/ \dd\cle thin then 
there is an ordinal\/ $\al<\ck$ and a\/ $\id11$ \lr\ order 
preserving map\/ $F:\stk{\bn}{\cle}\to\stk{2^\al}{\lexe}$ 
satisfying the following additional requirement$:$ 
if\/ $x,y\in\Xa$ then\/ ${x\napr y}\limp{F(x)\ne F(y)}\;;$ 

\itla{hmst2} 
if\/ $\Xa$ is\/ \dd\cle thin then\/ 
$\Xa$ is covered by the (countable) union of\/ 
all\/ $\id11$ \dd\cle chains\/ $C\sq\bn\;.$ 
\een
\ete

\lr\ order preserving maps $F$, satisfying the extra 
requirement of non-gluing of \dd\apr classes as in \ref{hmst1},  
were called \rit{linearization maps} in \cite{k:blin}. 

Any map $F$ as in \ref{hmst1} of the theorem sends 
any two \dd\cle in\-comparable reals $x,y\in \bn$ onto a 
\dd\lex comparable pair of $F(x),F(y)$, that is, either 
strictly $F(x)\lex F(y)$ or strictly $F(y)\lex F(x)$.
On the other hand, if the background set $\Xa$ is already 
a \dd\cle chain then $F$ has to be \rl\ order preserving  
too, that is, $x\cle y$ iff $F(x)\lexe F(y)$ for all 
$x,y\in\Xa$.

\bpf
[Claim \ref{mt1} of Theorem \ref{mt} modulo Theorem~\ref{hmst}] 
First of all, assume that the given Borel order $\cle$ is 
in fact $\id11$ and the given set $X=\Xa$ is $\is11$. 
The 
case of $\id11(p)$ and $\is11(p)$ with any fixed real parameter 
$p$ is accordingly reducible to a corresponding version of 
Theorem~\ref{hmst}.
 
Let, by Theorem~\ref{hmst}\ref{hmst2}, 
$\Xa\sq\bigcup_n C_n$, where each $C_n$ 
is a $\id11$ \dd\cle chain, and let $F$ and $\al$ 
be given by Theorem~\ref{hmst}\ref{hmst1}. 
To check that $\Xa$ is countably cofinal, it suffices to show 
that such is every set $X_n=\Xa\cap C_n$. 
But $X_n$ is a chain, so if it is {\ubf not} countably 
cofinal then there is a strictly \dd\cl increasing 
sequence $\sis{x_\al}{\al<\omi}$ of elements $x_\al\in X_n$. 
Then $\sis{F(x_\al)}{\al<\omi}$ is accordingly a strictly 
\dd\lex increasing sequence in $2^\al,$ which is impossible. 

Finally if there is a \dd\cle chain in $X$ of uncountable 
cofinality then a similar argument leads to such a chain in 
$\ang{2^\al;\lexe}$, with the same contradiction.\vom 

\epF{Theorem \ref{mt}\ref{mt1}}

Theorem~\ref{hmst} itself  
is an extension of two results in \cite{hms} 
(theorems 3.1 and 5.1). 
The latter directly correspond to the case of 
$\id11$ sets $\Xa$ in Theorem~\ref{hmst}. 
However the proof of Theorem~\ref{hmst} we 
manufactured in \cite{kanHMS} 
rather strictly follows the arguments in \cite{hms}. 
See also \cite{k:blin} in matters of the additional requirement in 
claim \ref{hmst1}, which also is presented in \cite{hms} 
implicitly.

\punk{Remarks and corollaries}
\las{sek1c} 

Claim~\ref{hmst1} of Theorem~\ref{hmst} can be strengthened 
as follows:
{\it
\ben
\Renup
\itla{lin1} 
if there is no continuous 1-1 \lr\ order preserving map\/ 
$F:{\stk{\dn}{\le_0}}\to\stk\Xa\cle$ 
such that\/ $a\not\Eo b$ implies that\/ $F(a)\yi F(b)$ 
are\/ \dd\cle incomparable, 
then there is an ordinal\/ $\al<\ck$ and a\/ $\id11$ \lr\ order 
preserving map\/ $F:\stk{\bn}{\cle}\to\stk{2^\al}{\lexe}$ 
satisfying the following additional requirement$:$ 
if\/ $x,y\in\Xa$ then\/ ${x\napr y}\limp{F(x)\ne F(y)}\;.$ 
\een
}
Here $\le_0$ is the PQO on $\dn$ defined so that $x\le_0y$ 
iff $x\Eo y$ and either $x=y$ or $x(k)<y(k)$, where $k$ is 
the largest number with $x(k)\ne y(k)$.\snos
{
$<_0$ orders each \dd\Eo class similarly to 
the (positive and negative) integers, except for the class 
$\ek{\om\ti\ans0}{\Eo}$ ordered as $\om$ and
the class $\ek{\om\ti\ans1}{\Eo}$ ordered the inverse of $\om$.} 
The ``if'' premice in \ref{lin1} is an immediate consequence 
of the \dd\cle thinness of $\Xa$ as in \ref{hmst1}, and hence 
\ref{lin1} really strenthens \ref{hmst1} of Theorem~\ref{hmst}.

Claim \ref{lin1} is an extension of Theorem 3 in
\cite{k:blin}; the latter corresponds to the case of $\id11$ 
sets $\Xa$. 

In the category of chains (rather than thin sets), 
the case of $\is11$ sets $\Xa$ in Theorem \ref{mt}\ref{mt1} 
is reducible to the case of $\id11$ 
sets simply because any $\is11$ chain $X$ can be covered by 
a $\id11$ chain $Y$. 
We find such a set $Y$ by means of the following two-step 
procedure.\snos
{See a different argument, based on 
a reflection principle in \cite[Corollary 1.5]{hms}.}
The set $C$ of all elements, 
\dd\cle comparable with every element $x \in X$, 
is $\ip11$,  
and $X \sq C$ (as $X$ is a chain).
By the Separation theorem, there is a $\id11$ set 
$B$ such that $X \sq B \sq C$. 
Now, the set $U$ of all elements in 
$B$, comparable with every element in $B$, is $\ip11$, 
and we have $X \sq B$. 
Once again, by Separation, there is a $\id11$ set $Y$ such 
that $X \sq Y \sq U$.
By construction, $U$ and $Y$ are chains, as required. 

Recall the following well-known earlier result in passing by, 
originally due to H.\,Friedman, as mentioned in 
\cite{hash}. 

\bcor
[of Theorem \ref{mt}\ref{mt1}] 
\lam{lea}
Every Borel LQO\/ $\le$ is countably cofinal, 
and moreover, there is no strictly increasing\/
\dd\omi sequences.\qed
\ecor 

The next immediate corollary says that maximal chains cannot 
be analytic provided they are not countably cofinal. 

\bcor
\lam{cora}
If\/ $\cle$ is a Borel PQO, and every countable set\/ 
$D\sq \dom{(\cle)}$ has a strict upper bound,
then there is no\/ $\fs11$ maximal\/ \dd\cle chains.\qed
\ecor 

\bcor
[Harrington and Shelah \cite{hash,shel}]
\lam{fn2}%
If\/ $\cle$ is a\/ $\fp11$ LQO on a Borel set  
then there is no strictly increasing\/ \ddd\omi chains in\/ $\cle$.
\ecor
\bpf   
The result was first obtained by a direct and rather 
complicated argument.
But fortunately there is a reduction to the Borel case.

Indeed let $x\cl y$ iff $y\not\cle x$, so in fact 
$R_0={\cl}$ is just the strict LQO associated with $\cle$.
As $R_0\sq(\cle)$, by Separation there is a Borel set $B_0$, 
$R_0\sq B_0\sq(\cle)$. 
Let $B_0'$ be the relation of \ddd{B_0}incomparability, and 
let $R_1$ be the PQO-hull of $B_0\cup B_0'$. 
Thus $R_1$ is a LQO and $R_0\sq B_0\sq R_1\sq (\cle)$. 

Once again, let $B_1$ is Borel set such that 
$R_1\sq B_1\sq(\cle)$. 
Define sets $B'_1$ and $R_2$ as above. 
And so on. 

Finally, after $\om$ steps, the union 
$R=\bigcup_nB_n=\bigcup_nR_n$ is a Borel LQO and 
$(\cl)\sq R\sq (\cle)$. 
Any strictly \ddd\cle increasing chain is 
strictly \ddd Rincreasing as well. 
It remains to apply Corollary~\ref{lea}.
\epf

\punk{Near-counterexamples for chains}
\las{sec2}

The following examples show that, 
even in the particular case of 
chains instead of thin orderings,  
Theorem \ref{mt}\ref{mt1} is not true any more 
for different extensions of the domain of 
$\fs11$ suborders of a Borel partial quasi-orders, 
such as $\fs11$ and $\fp11$ linear quasi-orders --- 
not necessarily suborders of Borel orderings, as well as 
$\fd12$ and $\fp11$ suborders of Borel orderings.
In each of these classes, a counterexample of cofinality 
$\omi$ will be defined.

\bex
[$\fs11$ LQO] 
\lam{ex1}
Consider a recursive coding of sets of rationals by reals. 
Let $Q_x$ be the set coded by a real $x$.
Let $X_\al$ be the set of all reals $x$ such that the 
maximal well-ordered initial segment
of $Q_x$ has the order type $\al$. 
We define 
$$
x\le y
\quad\text{iff}\quad
\sus\al\,\sus\ba\:(x\in X_\al\land y\in X_\ba\land\al\le\ba).
$$
Then $\le$ is a $\is11$ LQO on $\bn$ of cofinality $\om_1$.

Note that the associated strict order $x<y$,
iff $x\le y$ but not $y\le x$, 
is then more complicated than just $\fs11$,
therefore there is no contradiction in this example 
to the result mentioned in Remark~\ref{fn2}.\qed
\eex

\bex
[$\fp11$ LQO] 
\lam{ex2}
Let $D\sq\bn$ be the $\ip11$ set of codes of (countable) ordinals. 
Then
$$
x\le y
\quad\text{iff}\quad
x,y\in D\,\land\, |x|\le |y|
$$ 
is a $\ip11$ LQO of cofinality $\omi$.
Note that $\le$ is defined on a non-Borel $\ip11$ set $D$, 
and there is no $\fp11$ LQO of cofinality $\omi$ but 
defined on a Borel set --- by exactly the same argument as 
in Remark~\ref{fn2}.
\qed
\eex

\bex
[$\fp11$ LO] 
\lam{ex3}
To sharpen Example \ref{ex2}, define
$$
x\le y
\quad\text{iff}\quad
x,y\in D\;\land\;
\big( {|x|<|y|}\,\lor\, 
{(|x|=|y|\land x\lex y)}\big);
$$ 
this is a $\ip11$ LO of cofinality $\omi$.\qed
\eex

\bex
[$\fd12$ suborders] 
\lam{ex4}
Let $\le$ be the eventual domination order on $\om^\om$. 
Assuming the axiom of constructibility $\rV=\rL$, one 
can define a 
strictly $\le$-increasing $\id12$ $\omi$-sequence 
$\sis{x_\al}{\al<\omi}$ in $\om^\om$.\qed 
\eex

\bex
[$\fp11$ suborders] 
\lam{ex5}
Define a PQO $\le$ on $(\om\bez\ans0)^\om$ so that
$$
x\le y
\quad\text{ iff }\quad
\text{either}\quad x=y \quad\text{or}\quad
\lim_{n\to\iy}\:y(n)\,{/\hspace{-1ex}/}\,x(n)=\iy 
$$
(the ``or'' option defines the associated strict order $<$).
Assuming the axiom of constructibility $\rV=\rL$, define a 
strictly increasing $\id12$ $\omi$-sequence 
$\sis{x_\al}{\al<\omi}$ in $\om^\om$. 
By the Novikov -- Kondo -- Addison $\ip11$ Uniformization 
theorem, there is a $\ip11$ set 
$\sis{\ang{x_\al,y_\al}}{\al<\omi}\sq\bn\ti\dn$. 
%
%
%
Let $z_\al(n)=3^{x_\al(n)}\cdot 2^{y_\al(n)}$, $\kaz n$.
Then the \ddd\omi sequence $\sis{z_\al}{\al<\omi}$ is $\ip11$ 
and strictly increasing: indeed, factors of the form 
$2^{y_\al(n)}$ are equal 1 or 2 whenever $\al\in\dn$.\qed
\eex

\punk{Definable thin suborders in the Solovay model}
\las{sbSM}

Here {\ubf we prove Theorem \ref{mt}\ref{mt2}}. 
Arguing in the Solovay model 
(a model of \ZFC\ defined in \cite{solmod}, in which all \ROD\ 
sets of reals are Lebesgue measurable), 
we assume that $\cle$ is a Borel PQO on a Borel set 
$D\sq\bn$,  
$X\sq D$ is a $\ROD$ (real-ordinal definable) set, 
and the set $X$ is a \dd\cle thin. 

Let $\rho<\omi$ be such that $\cle$ is a relation in $\fs0\rho$. 

Prove that 
the restricted ordering $\ang{X;{\cle}}$ 
\rit{is countably cofinal}, \ie, contains a countable 
cofinal subset (not necessarily a chain). 

It is known that in the Solovay model any \ROD\ set in $\bn$ 
is a union of a \ROD\ \ddd\omi sequence of analytic sets. 
Thus there is a \ddd\sq increasing \ROD\ sequence 
$\sis{X_\al}{\al<\omi}$ of $\fs11$ sets $X_\al$, such that 
$X=\bigcup_{\al<\omi}X_\al$. 
Let $r\in\bn$ be a real parameter such that in fact  the 
sequence $\sis{X_\al}{\al<\omi}$ is $\od(r)$.

As the sets $X_\al$ are countably \dd\cle cofinal by claim 
\ref{mt1} of Theorem~\ref{mt}, 
it suffices to prove that one of $X_\al$ is cofinal in $X$. 
 
{\it Suppose otherwise\/}. 
Then the sets 
${D_\al=\ens{z\in D}{\sus x\in X_\al\,(z\cle x )}}$  
contain $\ali$ different sets and form an $\od(r)$ sequence. 
We claim that 
every set $D_\al$ belongs to the same class  
$\fs0\rho$ as the given Borel order  $\cle$.
Indeed let 
$\ens{x_n}{n\in\om}$ 
be any countable cofinal set in $X_\al$. 
Then the set $D_\al=\ens{z\in D}{\sus n\,(z\cle x_n)}$ 
is $\fs0\rho$ by obvious reasons.%

We conclude that the Borel class $\fs0\rho$ contains $\ali$ 
pairwise different sets in $\od(r)$ for one and the same 
$r\in\bn.$ 
But this contradicts to a well-known result of 
Stern \cite{stern}.

\qeDD{Theorem \ref{mt}\ref{mt2}}

\punk{$\fs12$ thin suborders of Borel PQOs}
\las{sbS12}

Here {\ubf we prove Theorem \ref{mt}\ref{mt3}}. 
Assume that $\cle$ is a Borel PQO on a Borel set $D\sq\bn$,  
$X\sq D$ is a $\fs12$ set, 
and\/ $X$ is \dd\cle thin.  
We also assume that $\omi^{\rL[r]}<\omi$ for every real $r$. 

Prove that the ordering $\ang{X;{\cle}}$ 
\rit{is countably cofinal}.

Pick a real $r$ such that $X$ is $\is12(r)$ and 
$\cle$ is $\id11(r)$.
To prepare for an absoluteness argument, 
fix canonical formulas,
$$
\vpi(\cdot,\cdot)\;\; \text{of type}\;\;\is12\,,\quad 
\sg(\cdot,\cdot,\cdot)\;\; \text{of type}\;\;\is11\,,\quad 
\pi(\cdot,\cdot,\cdot)\;\; \text{of type}\;\;\ip11\,, 
$$
which define $X$ and $\cle$ in the set universe $\rV$, 
so that it is true
in $\rV$ that
$$
x\cle y\leqv \sg(r,x,y)\leqv \pi(r,x,y)
\quad\text{and}\quad
x\in X\leqv\vpi(r,x)\,.
$$
for all $x,y\in\bn.$
We let
$X_{\vpi}=\ens{x\in\bn}{\vpi(r,x)}$
and
$$
x\le_{\sg\pi} y
\quad\eqv\quad \sg(r,x,y)
\quad\eqv\quad \pi(r,x,y)
$$
so that $X_{\vpi}=X$ and $\le_{\sg\pi}$ is $\cle$ in $\rV$, but
$X_{\vpi}$ and $\le_{\sg\pi}$ can be defined 
in any transitive universe containing $r$ 
and containing all ordinals 
(to preserve the equivalence of formulas $\sg$ and $\pi$).

Let $\wo$ be the canonical $\ip11$ set of codes of (countable) 
ordinals, and for $w\in\wo$ let $|w|<\omi$ be the ordinal 
coded by $w$.

Let $X_\vpi=\bigcup_{\al<\omi}X_\vpi(\al)$ be a canonical 
representation of $X_\vpi$ as an increasing union of $\fs11$ sets.
Thus to define $X_\vpi(\al)$ fix a $\ip11(r)$ set 
$P\sq(\bn){}^2$ such that $X=\ens{x}{\sus y\,P(x,y)}$, 
fix a canonical $\ip11(r)$ norm 
$f:P\to\omi$, and let 
$$
P_\al=\ens{\ang{x,y}}{f(x,y)<\al}
\quad\text{and}\quad
X_\vpi(\al)=\ens{x}{\sus y\,(\ang{x,y}\in P_\al)}\,.
$$

In our assumptions, the ordinal $\Omega=\omi$ is inaccessible 
in $\rL[r]$. 
Let $\cP=\text{Coll}({{<}\,\Om},\om)\in \rL[r]$ 
be the corresponding Levy collapse forcing.  
Consider a \ddd\cP generic extension 
$\rV[G]$ of the universe.  
Then $\rL[r][G]$ is a Solovay-model generic extension 
of $\rL[r]$.  
The plan is to compare the models $\rV$ and $\rL[r][G]$. 
Note that $\rL[r]$ is their common part, 
$\rV[G]$ is their common extension, and the 
three models have the same cardinal 
$\omi^\rV=\omi^{\rL[r][G]}=\omi^{\rV[G]}=\Om>\omi^{\rL[r]}$. 

\ble
\lam{parta}
It is true both in\/ $\rV[G]$ and\/ $\rL[r][G]$ that if\/ 
$\al<\Om$ then the set\/ $X_\vpi(\al)$ is\/ 
\ddd{\le_{\sg\pi}}thin.
\ele

This key absoluteness lemma has no analogies in a 
simpler case of chains (instead of thin sets) 
earlier considered in \cite{cocof}.
Here we even don't claim the absolutenes of the thinness 
poperty of the whole set  
$X_\vpi=\bigcup_{\al<\Om}X_\vpi(\al)$!

\bpf[Lemma]
Note that the thinness of $X_\vpi(\al)$ is a $\ip13$ 
statement with parameters $r$ and any real which codes $\al$. 
This makes the step $\rV[G]\to\rL[r][G]$ trivial 
by Shoenfield, and 
allows to concentrate on $\rV[G]$.

Suppose towards the contrary that there is a perfect tree 
$T\in\rV[G]$, $T\sq\nse,$ such that the perfect set 
$[T]=\ens{x\in\bn}{\kaz n\,(x\res n\in T)}$ satisfies 
\ben
\aenu
\itla{par1}
$[T]\sq X_\vpi(\al)$ and $[T]$ is a \ddd{\le_{\sg\pi}}antichain 
\een
in $\rV[G]$. 
There exist an ordinal $\ga<\Om$ and a 
\ddd{\col\ga}generic map $F\in\rV[G]$ 
such that already $T\in\rV[F]$, so that $T=t[F]$, where 
$t\in\rV$, $t\sq {\col\ga}\ti\nse$ is a \ddd{\col\ga}name.

Note that \ref{par1} is still true in $\rV[F]\sq\rV[G]$ by 
Shoenfield, moreover, \ref{par1} is true in $\rL[z][F]$, 
where a real $x\in\bn\cap\rV$ codes all of $\al,\ga,r,t$.
Therefore there is a condition $s\su F$ 
(a finite string of ordinals $\xi<\ga$) 
which \ddd{\col\ga}forces \ref{par1}
(with $T$ replaced by the name $t$)
over $\rL[z]$. 

Now, by the assumptions of Theorem \ref{mt}\ref{mt3}, 
there is a map $F'\in\rV$, still \ddd{\col\ga}generic
over $\rL[z]$ and satisfying $s\su F'$. 
Then the tree $T'=t[F']$ belongs to the model 
$\rL[z][F']\sq\rV$ and satisfies \ref{par1} 
(in the place of $T$) in $\rL[z][F']$, hence, in $\rV$
as well by Shoenfield. 
But this contradicts to the choice of $X=X_\vpi$.
\epF{Lemma}

We continue the proof of Theorem \ref{mt}\ref{mt3}. 
It follows from the lemma that all orderings 
$\stk{X_{\vpi}(\al)}{\le_{\sg\pi}}$, $\al<\Om$, 
are countably cofinal in $\rL[r][G]$ by 
Theorem~\ref{mt}\ref{mt1}. 
However $\rL[r][G]$ is a Solovay-model type extension 
of $\rL[r]$. 
Therefore (see the argument in Section~\ref{sbSM}) 
it is true in $\rL[r][G]$ that the whole ordering  
$\stk{X_{\vpi}}{\le_{\sg\pi}}$ is countably cofinal, 
hence there is an ordinal $\al<\Om=\omi^{\rL[r][G]}$ 
such that the sentence 
\ben
\aenu
\atc
\itla{f1}
the subset $X_\vpi(\al)$ is 
\ddd{\le_{\sg\pi}}cofinal in the whole set $X_\vpi$ 
\een
is true in $\rL[r][G]$.
However \ref{f1} can be expressed by a $\ip12$ formula with 
$r$ and an arbitrary code $w\in\wo\cap\rL[r][G]$ such that 
$|w|=\al$ --- as the only parameters. 
It follows, by Shoenfield, that \ref{f1} 
is true in $\rV[G]$ as well.

Then by exactly the same absoluteness argument  
\ref{f1} is true in $\rV$, too.
Thus it is true in $\rV$ that $X_\vpi(\al)$, 
a $\fs11$ set, is cofinal in the whole set $X= X_\vpi$. 
But $X_\vpi(\al)$ is countably cofinal by 
Theorem \ref{mt}\ref{mt1}.\vom

\qeDD{Theorem \ref{mt}\ref{mt3}}

\punk{The Solovay model and $\od$ forcing}
\las{smod}

Here we begin the proof of Theorem~\ref{t:sm}.
We emulate the proof of Theorem~5.1 in \cite{hms} and 
a similar proof of Theorem~\ref{hmst}\ref{hmst2} above 
(given in \cite{kanHMS}), changing 
the Gandy -- Harrington forcing $\gh$ with the $\od$ forcing $\OD$.
There is no direct analogy between the two forcing notions, so 
we'll both enjoy some simplifications and suffer from some 
complications.

We start with a brief review of the Solovay model.
Let $\Om$ be an ordinal. 
Let \osm\ be the following hypothesis:
\bde
\item[\rm\osm:]
$\Om=\omi$, 
$\Om$ is strongly inaccessible in $\rL$, the 
constructible universe,   and
the whole universe $\rV$ is 
a generic extension of $\rL$ via the Levy collapse forcing 
$\coll\Om$, as in \cite{solmod}.
\ede
Assuming \osm, let $\OD$ be the set of all {\ubf non-empty} 
$\od$ sets $Y\sq \bn$.
We consider $\OD$ as a forcing notion 
(smaller sets are stronger).
A set $D\sq \OD$ is:
\bit
\item[$-$]\rit{dense}, \ 
iff for every $Y\in\OD$ there exists $Z\in D$, $Z\sq Y$;

\item[$-$]\rit{open dense}, \ 
iff in addition we have $Y\in D\imp X\in D$ whenever 
sets $Y\sq X$ belong to $\OD$;
\eit   
A set $G\sq \OD$ is {\pge}, \ iff \ 
1) 
if $X,Y\in G$ then there is a set $Z\in G$, $Z\sq X\cap Y$,   
\hspace*{0.3ex}and \ 
2) 
if $D\sq\OD$ is $\od$ and dense then $G\cap D\ne\pu$.

Given an $\od$ \eqr\ $\rE$ on $\bn,$ a 
\rit{reduced product} forcing 
notion $\spe$ consists of all sets of the form 
$X\ti Y,$ where $X\yi Y\in\OD$ and $\eke X\cap\eke Y\ne\pu$.
For instance $X\ti X$ belongs to $\spe$ whenever $X\in\OD$. 
The notions of sets dense and open dense in $\spe$, and 
\ege\ sets are similar to the case of $\OD$

A condition $X\ti Y$ in $\spe$ is \rit{saturated} iff 
$\eke X=\eke Y$.

\ble
\lam{sm6}
If\/ $X\ti Y$ is a condition in $\spe$ then there is a 
stronger saturated subcondition\/ $X'\ti Y'$ in $\spe$. 
\ele
\bpf
Let $X'=X\cap\eke Y$ and $Y'=Y\cap\eke X$. 
\epf

\bpro
[lemmas 14, 16 in \cite{ksol}] 
\lam{genx}
Assume \osm. 

If a set\/ $G\sq\OD$ is\/ \pge\ then the 
intersection\/ $\bigcap G=\ans{x[G]}$ consists of a 
single real\/ $x[G]$, called\/ \pge\ ---  
its name will be\/ $\dox$.

Given an $\od$ \eqr\ $\rE$ on $\bn,$ 
if\/ $G\sq\spe$ is\/ \ege\ then the 
intersection\/ $\bigcap G=\ans{\ang{\xle[G],\xri[G]}}$ 
consists of a single pair of reals\/ $\xle[G]\yi\xri[G]$, 
called an\/ \ege\ pair --- 
their names will be\/ $\doxl\yi\doxr$; 
either of\/ $\xle[G]\yi\xri[G]$ is separately\/ \pge.
\qed
\epro

As the set $\OD$ is definitely uncountable, 
the existence of \pge\ sets does not 
immediately follow from \osm\ by a cardinality argument. 
Yet fortunately $\OD$ is \rit{locally countable}, in a sense. 

\bdf
[assuming \osm]
\lam{odik}
A set $X\in\od$ is \rit{\odk} if the set 
$\pwod X=\pws X\cap\od$ of all $\od$ subsets of $X$ 
is at most countable. 
\edf

For instance, assuming \osm, the set $X=\bn\cap\od=\bn\cap\rL$ 
of all $\od$ reals is \odk. 
Indeed  
$\pwod X = \pws X\cap \rL$, and hence 
$\pwod X$ admits an $\od$  
bijection onto the ordinal $\om_2^{\rL}<\omi=\Om$.

\ble
[assuming \osm]
\lam{2co}
If a set $X\in\od$ is \rit{\odk} then the set\/ $\pwod X$ 
is\/ \odk either.
\ele
\bpf
There is an ordinal $\la<\omi=\Om$ and an $\od$ bijection 
$b:\la\onto\pwod X$. 
Any $\od$ set $Y\sq\la$ belongs to $\bL$, hence, the $\od$ 
power set $\pwod\la=\pws\la\cap\bL$ belongs 
to $\bL$ and $\card(\pwod\la)\le\la^+<\Om$ in $\bL$. 
We conclude that $\pwod\la$ is countable. 
It follows that $\pwod{\pwod X}$ is countable, as required.
\epf


\ble
[assuming \osm]
\lam{L*}
If\/ $\la<\Om$ then the set\/ $\koh_\la$ of all 
elements\/ $f\in\la^\om$, \dd{\col\la}generic over\/ $\rL$, 
is\/ \odk. 
\ele
\bpf
If $Y\sq\koh_\la$ is \od\ and $x\in Y$ then 
``$\dotx\in \dotY$'' is \dd{\col\la}forced over $\rL$. 
It follows that there is a set 
$S\sq\la\lom=\col\la\yt S\in\rL$, such 
that $Y=\koh_\la\cap\bigcup_{t\in S}\ibn t$, where 
$\ibn t=\ens{x\in\la\lom}{t\su x}$, 
a Baire interval in $\la\lom.$ 
But the collection of all such sets $S$ belongs to $\rL$ 
and has cardinality 
$\la^+$ in $\rL$, hence, is countable under \osm.
\epf

Let $\odi$ be the set of all \odk\ sets $X\in\OD$.
We also define 
$$
\spei=\ens{X\ti Y\in\spe}{X,Y\in \odi}.
$$ 

\ble
[assuming \osm]
\lam{den}
The set\/ $\odi$ is dense in\/ $\OD$, that is, 
if\/ $X\in\OD$ then there is a condition\/ $Y\in\odi$ 
such that\/ $Y\sq X$. 

If\/ $\rE$ is an\/ $\od$ \eqr\ on $\bn$ 
then the set\/ $\spei$ is dense in\/ $\spe$ 
and any\/ $X\ti Y$ in\/ $\spei$ is\/ \odk.
\ele
\bpf
Let $X\in\OD$. 
Consider any $x\in X$. 
It follows from \osm\ that there is 
an ordinal $\la<\omi=\Om$, 
an element $f\in\koh_\la$, 
and an $\od$ map $H:\la^\om\to\bn$, 
such that $x=H(f)$. 
The set $P=\ens{f'\in\koh_\la}{H(f')\in X}$ is then \od\ 
and non-empty (contains $f$), and hence so is its image 
$Y=\ens{H(f')}{f'\in P}\sq X$ (contains $x$). 
Finally, $Y\in\odi$ by Lemma~\ref{L*}. 

To prove the second claim, let $X\ti Y$ be a condition in 
$\spe.$ 
By Lemma~\ref{sm6} there is a stronger saturated 
subcondition $X'\ti Y'\sq X\ti Y$.
By the first part of the lemma, let $X''\sq X'$ be a 
condition in $\odi$, and $Y''=Y'\cap\eke{X''}$.
Similarly, let $Y'''\sq Y''$ be a 
condition in $\odi$, and
$X'''=X''\cap\eke{Y'''}$.
Then $X'''\ti Y'''$ belongs to $\spei$.
\epf

\bcor
[assuming \osm]
\lam{egen}
If\/ $X\in\OD$ then there exists a \pge\ set\/ $G\sq\OD$ 
containing\/ $X$. 
If\/ $X\ti Y$ is a condition in $\spe$ then there exists 
a \ege\ set\/ $G\sq\spe$ containing\/ $X\ti Y$. 
\ecor
\bpf
By Lemma~\ref{den}, assume that $X\in\odi$. 
Then the set $\OD_{\sq X}$ of stronger conditions 
contains only countably many \od\ subsets by 
Lemma~\ref{2co}.
\epf

\punk{The $\od$ forcing relation}
\las{odfr}
        
The forcing notion $\OD$ will play the same role below as 
the Gandy -- Harrington forcing in \cite{hms,kanHMS}.
There is a notable technical difference: under \osm, 
\dodg\ sets exist in the ground Solovay-model universe 
 by Corollary~\ref{egen}. 
Another notable difference is connected with the forcing 
relation. 

\bdf
[assuming \osm]
\lam{frd}
Let $\vpi(x)$ be an \dd\Ord\rit{formula}, that is, a 
formula with ordinals as parameters.

A  condition $X\in \OD$ is said to \rit{\pfo} $\vpi(\dox)$
iff $\vpi(x)$ is true 
(in the Solovay-model set universe considered)
for any \pge\ real $x$. 

If $\rE$ is an\/ $\od$ \eqr\ on $\bn$ then 
a condition $X\ti Y$ in $\spe$ is said to 
\rit{\efo} $\vpi(\doxl,\doxr)$
iff $\vpi(x,y)$ is true 
for any \ege\ pair $\ang{x,y}$.
\edf

\ble
[assuming \osm]
\lam{frl}
Given an\/ \dd\Ord formula\/ $\vpi(x)$ and a\/ 
\pge\ real\/ $x$, if\/ $\vpi(x)$ is true\/ 
{\rm (in the Solovay-model set universe considered)} 
then there is a condition $X\in \OD$ containing\/ $x$, 
which\/ \pfo s\/ $\vpi(\dox)$.

Let\/ $\rE$ be an\/ $\od$ \eqr\ on $\bn.$ 
Given an\/ \dd\Ord formula\/ $\vpi(x,y)$ and a\/ 
\ege\ pair\/ $\ang{x,y}$, if\/ $\vpi(x,y)$ is true 
then there is a condition in\/ $\spe$ 
containing\/ $\ang{x,y}$, 
which\/ \efo s\/ $\vpi(\doxl,\doxr)$.
\ele
\bpf
To prove the first claim, put $X=\ens{x'\in\bn}{\vpi(x')}$. 
But this argument does not work for $\spe$. 
To fix the problem, we propose a longer argument which 
equally works in both cases --- but we present it in the 
case of $\OD$ which is slightly simpler.

Formally the forcing notion $\OD$ does not belong to $\rL$. 
But it is order-isomorphic to a certain forcing notion 
$P\in\rL$, namely, the set $P$ of \rit{codes}\snos
{A code of an $\od$ set $X$ is a finite sequence of logical 
symbols and ordinals which correspond to a definition in the 
form $X=\ens{x\in\rV_\al}{\rV_\al\models\vpi(x)}$.} 
of $\od$ sets in $\OD$.
The order between the codes in $P$, 
which reflects the relation $\sq$ 
between the $\od$ sets themselves, is expressible in $\rL$, 
too. 
Furthermore dense $\od$ sets in $\OD$ correspond 
to dense sets in the coded forcing $P$ in $\rL$.

Now, let $x$ be \pge\ and $\vpi(x)$ be true. 
It is a known property of the Solovay model that there is 
another \dd\Ord formula\/ $\psi(x)$ such that 
$\vpi(x)$ iff $\rL[x]\models\psi(x)$. 
Let $g\sq P$ be the set of all codes of conditions $X\in\OD$ 
such that $x\in X$. 
Then $g$ is \dd Pgeneric over $\rL$ by the choice of $x$, 
and $x$ is the corresponding generic object, 
hence there is a condition $p\in g$ which \dd Pforces 
$\psi(\dox)$ over $\rL$. 
Let $X\in\OD$ be the \od\ set coded by $p$, so $x\in X$. 
To prove that $X$ \odf s $\vpi(\dox)$, let $x'\in X$ be a 
\pge\ real.
Let $g'\sq P$ be the \dd Pgeneric set of all codes of 
conditions $Y\in\OD$ such that $x'\in Y$. 
Then $p\in g'$, hence
$\psi(x')$ holds in $\rL[x']$, by the choice 
of $p$. 
Then $\vpi(x')$ holds 
(in the Solovay-model set universe) by the choice of $\psi$, 
as required.
\epf

\bcor
[assuming \osm]
\lam{frc}
Given an\/ \dd\Ord formula\/ $\vpi(x)$, if\/ $X\in\OD$ 
does not\/ \pfo\/ $\vpi(\dox)$ then there is a condition 
$Y\in \OD\yd Y\sq X$, 
which\/ \pfo s\/ $\neg\:\vpi(\dox)$.
The same for\/ $\spe$.\qed 
\ecor

\punk{Adding a perfect antichain}
\las{fpis}

The next result will be pretty important.

\ble
[assuming \osm]
\lam{29sm}
Assume that\/ $\cle$ is an\/ $\od$ PQO on\/ $\bn,$
and\/ $\rE_A$ is an\/ $\od$ \eqr\ on\/ $\bn$ 
for any\/ $A\in\OD$, such that if\/ $A\sq B$ then\/ 
$x\rE_A y$ implies\/ $x\rE_B y$.

Suppose that\/ $\Xa\in\OD$, and if\/ $B\in\OD\yt B\sq \Xa$ 
then\/ $B\ti B$ does  
{\bfit not}\/ \efa B\ 
that\/ $\doxl,\doxr$ are\/ \dd\cle comparable. 
Then\/ $\Xa$ is not\/ \dd\cle thin. 
\ele
\bpf
[follows 2.9 in \cite{hms}]
Let $T$ be the set of all finite trees $t\sq2\lom.$
If $t\in T$ then let $M(t)$ be the set of all \dd\su maximal 
elements of $t$. 

Let $\Phi$ be the set of systems $\vpi=\sis{X_u}{u\in t}$ of sets 
$X_u\in\odi$, such that $t\in T$ and the 
following conditions \ref{291} -- \ref{294} are satisfied:
\ben
\renu
\itla{291}
$X_\La \sq X^\ast$ (where $\La$ is the empty string);

\itla{292}
if $u\su v\in t$ then $X_v\sq X_u$;

\itla{293} 
if $u\we0$ and $u\we 1$ belong to $t$ then 
$X_{u\we 0}\ti X_{u\we1}$ belongs to $\spai{X_u}$ 
and \efa{X_u}s that 
$\doxl$ is \dd\cle incomparable to $\doxr$;

\itla{294} 
compatibility: 
there is a sequence $\sis{x_u}{u\in M(t)}$ of points $x_u\in X_u$ 
such that if $u,v\in M(t)$ then $x_u\rE_{X_{u\land v}}x_v$, 
where $u\land v$ is the largest string $w\in2\lom$ such that 
$w\su u$ and $w\su v$ --- it easily follows that then 
$X_u\ti X_v$ is a condition in $\spa{X_{u\land v}}$.
\een 
Say that a system $\sis{X_u}{u\in t}\in\Phi$ is 
\rit{saturated} if in addition
\ben
\renu
\atc\atc\atc\atc  
\itla{295}
for any $v\in M(t)$ and $x\in X_v$ there is a sequence 
$\sis{x_u}{u\in M(t)}$ as in \ref{294}, such that $x_v=x$. 
\een

Say that a system $\sis{X'_u}{u\in t'}\in\Phi$: 
1) 
\rit{weakly extends} another system 
$\vpi=\sis{X_u}{u\in t}$ if $t\sq t'$,
$X_u=X'_u$ for all $u\in t\bez M(t)$, and  
$X'_u\sq X_u$ for all $u\in M(t)$; and 
2) 
\rit{properly extends} $\sis{X_u}{u\in t}$ if $t\sq t'$ 
and $X_u=X'_u$ for all $u\in t$.
Thus a weak extension not just adds new sets to a given 
system $\vpi$ but also shrinks old sets of the top layer 
$\vpi=\sis{X_u}{u\in M(t)}$ of $\vpi$.

\bcl
\lam{29d}
For any system\/ $\vpi=\sis{X_u}{u\in t}\in\Phi$ there is 
a {\bfit saturated\/} system\/ $\sis{X'_u}{u\in t}$ 
in $\Phi$ (with the same domain\/ $t$) which 
weakly extends\/ $\vpi$.
\ecl
\bpf 
If $u\in M(t)$ then simply let $X'_u$ be the set of all points 
$x\in X_u$ such that $x=x_u$ for some sequence 
$\sis{x_u}{u\in M(t)}$ as in \ref{294}.
\epF{claim}

\bcl
\lam{29c}
For any saturated system\/ 
$\vpi=\sis{X_u}{u\in t}\in\Phi$, 
if\/ $u\in M(t)$ then there are sets $X_{u\we0}\yi X_{u\we1}$ 
such that the system $\vpi$ extended 
by those sets still belongs to $\Phi$ and properly extends\/ 
$\vpi$.
\ecl
\bpf 
As $X_u\in\OD$ and $X_u\sq X^\ast$, the condition 
$X_u\ti X_u$ does\/ 
{\ubf not} \efa{X_u}\ 
that\/ $\doxl,\doxr$ are\/ \dd\cle comparable. 
By Corollary~\ref{frc}, 
pick a stronger condition $U\ti V\sq X_u\ti X_u$ in $\spa{X_u}$ 
which \efa{X_u}s 
that $\doxl,\doxr$ are \dd\cle incomparable.
By Lemmas \ref{den} and \ref{sm6} 
we may assume that $U\ti V$ belongs to $\spai{X_u}$ and is 
\dd{\rE_{X_u}}saturated, so that 
$\ek{U}{\rE_{X_u}}=\ek{V}{\rE_{X_u}}$.
We assert that the sets $X_{u\we0}=U$ and $X_{u\we1}=V$ prove 
the claim. 
It's enough to check \ref{294} for the extended system. 

Fix any $x\in X_{u\we0}=U$.
Then $x\in X_{u}$, hence, as the given system is saturated,
there is a sequence $\sis{x_v}{v\in M(t)}$ of points 
$x_v\in X_v$ as in \ref{294}, such that $x_u=x$. 
On the other hand, 
as $\ek{U}{\rE_{X_u}}=\ek{V}{\rE_{X_u}}$, there is a point 
$y\in V=X_{u\we1}$ such that $x\rE_{X_u}y$.
Put $x_{u\we0}=x$ and $x_{u\we1}=y$.
\epF{claim}

If $\rE$ is an $\od$ \eqr\ and 
${X\ti Y}\in {\spei}$ 
then the set $\cD(\rE,X,Y)$ of all 
sets, \rit{open dense in\/ $\spe$ below\/ $X\ti Y$}\snos
{\label{restf}%
That is, open dense subsets of the restricted forcing
$(\spe)_{\sq X\ti Y}=\ens{X'\ti Y'\in\spe}{X'\sq X\land Y'\sq Y}$.},
\imar{fnote restf}
is countable  by Lemma~\ref{den}; 
fix an enumeration 
$\cD(\rE,X,Y)=\ens{D_n(\rE,X,Y)}{n\in\om}$ such that 
$D_n(\rE,X,Y)\sq D_m(\rE,X,Y)$ whenever $m<n$. 

\bcl
\lam{29b}
Let\/ $n\in\om$ and\/ $\vpi=\sis{X_u}{u\in 2^{\le n}}\in\Phi$. 
Then there is a system\/ 
$\vpi'=\sis{X'_u}{u\in 2^{\le n+1}}\in\Phi$ which weakly 
extends\/ $\vpi$ and satisfies the following additional 
genericity requirement$:$ 
\ben
\fenu
\itla{**}
if strings\/ $u\ne v$ belong\/ to $2^{n+1}$ and\/ $w=u\land v$ 
{\rm(defined as in \ref{294})}  
then\/ 
condition\/ 
$X_{u}\ti X_{v}$ 
belongs to\/ $D_{n}(\rE_{X_w},X_{w\we0},X_{w\we1})$.
\een
\ecl
\bpf 
We first extend $\vpi$ by one layer of sets
$X'_{u\we i}$, $u\in 2^{n}$ and $i=0,1$, obtained by
consecutive $2^n$ splitting operations as in Claim \ref{29c},
followed by the saturating reduction as in Claim \ref{29d}. 
This way we get a saturated system  
$\eta=\sis{Y_u}{u\in 2^{\le n+1}}\in\Phi$ 
which weakly extends\/ $\vpi$. 

To fulfill \ref{**}, let us shrink 
the sets in the top layer $\sis{Y_u}{u\in 2^{n+1}}$ 
of $\eta$. 

Consider any pair of strings $u\ne v$ in $2^{n+1}.$
Let $w=u\land u$, so that $k=\dom w<n$, $w\su u$, $w\su v$, 
and $u(k)\ne v(k)$; let, say, $u(k)=0$, $v(k)=1$.
Condition $Y_{w\we0}\ti Y_{w\we1}$ belongs to 
$\spai{Y_w}$  by \ref{293} while $Y_{u}\ti Y_{v}$ belongs 
to $\spa{Y_w}$  by \ref{294} and satisfies 
$Y_{u}\sq Y_{w\we0}$ and $Y_{v}\sq Y_{w\we1}$ by \ref{292}.
By the density, there is a subcondition 
$Z_u\ti Z_v\sq Y_{u}\ti Y_{v}$ in 
$D_{n}(\rE_{Y_w},Y_{w\we0},Y_{w\we1})$; 
in particular, $Z_u\ti Z_v$ 
still belongs to $\spa{Y_w}$. 
In addition to $Z_u$ and $Z_v$, we 
let $Z_s=Y_s$ for any $s\in2^{n+1}\bez\ans{u,v}$. 
Then $\psi=\sis{Z_s}{s\in2^{\le n+1}}$ is still a system in 
$\Phi$. 
By Claim~\ref{29d}, there is a saturated system 
$\psi'=\sis{Z'_s}{s\in2^{\le n+1}}\in\Phi$ such that 
$Z'_s=Z_s=Y_s$ for all $u\in2^{\le n}$, and 
$Z'_s\sq Z_s$ for all $s\in2^{n+1}$. 
Then $Z'_u\sq Z_u$ and $Z'_v\sq Z_v$ --- so that 
$Z'_u\ti Z'_v\in D_{n}(\rE_{Y_w},Y_{w\we0},Y_{w\we1})$.

Iterating this shrinking construction $2^n(2^n-1)$ times 
(the number of pairs $s\ne t$ in $2^n$), 
we get a required system $\vpi'$. 
\epF{claim}

Claim~\ref{29b} allows to define, by induction, 
sets $X_u\sq X'_u\sq\Xa$ in $\odi$ ($u\in\bse$) and systems 
$\vpi_n=\sis{X_u}{u\in2^{<n}}\cup\sis{X'_u}{u\in2^{n}}$, 
such that, 
for any $n$: 
\ben
\aenu
\itla{vpi1}\msur
$\vpi_n$ is a saturated system in $\Phi$, 
weakly extended by $\vpi_{n+1}$, and 

\itla{vpi2}\msur
condition \ref{**} of Claim~\ref{29b} 
holds.
\een
Show that this leads to a required perfect set.

Suppose that $a\ne b$ are reals in $\dn,$ and $w=a\land b$, 
so that $w\su a$, $w\su b$, and $a(k)\ne b(k)$, where 
$k=\dom w$; let, say, $a(k)=0$, $b(k)=1$. 
Then the sequence of sets $X_{a\res m}\ti X_{b\res m}$, 
$m>k$, is \ega{X_w}\ by \ref{vpi2}, so that the intersection 
$\bigcap_{m>k}(X_{a\res m}\ti X_{b\res m})$ consists of a 
single pair of reals $\ang{x_a,x_b}$ 
by Proposition \ref{genx}.
Moreover, $x_a\yi x_b$ are \dd\cle incomparable by \ref{293}.
Finally it easily follows from \ref{vpi2} thet the diameters 
of sets $X_n$ uniformly tend to $0$ with $n\to\iy$, and hence 
the map $a\longmapsto x_a$ is continuous.
Thus $P=\ens{x_a}{a\in\dn}$ is a perfect 
\dd\cle antichain in $\Xa$.
\epf

\punk{Compression lemma}
\las{CL}

Let $\vT=\Om^+$; 
the cardinal successor of $\Om$ in both $\bL$, the ground model, 
and its \dd{\coll\Om}generic extension postulated by \osm\ to 
be the set universe; in the latter, $\Om=\omi$ and $\vT=\om_2$.

\ble
[compression lemma]
\lam{apal31}
Assume that\/ $\Om\le\vt\le\vT$ and\/ $X\sq2^\vT$ is the image of\/ 
$\bn$ via an\/ $\od$ map. 
Then there is an\/ $\od$ antichain\/ $A(X)\sq2\Lom$ and an\/ 
$\od$ isomorphism\/ $f:\stk{X}{\lexe}\onto \stk{A(X)}{\lexe}$.
\ele

Note that any antichain $A\sq2\Lom$ is linearly ordered by $\lexe$!

\bpf
If $\vt=\vT$ then, 
as $\card X\le\card{\bn}=\Om$, there is an ordinal $\vt<\vT$ 
such that $x\res\vt\ne y\res\vt$ whenever $x\ne y$ belong to $X$ 
--- this reduces the case $\vt=\vT$ to the case 
$\Om\le\vt<\vT$. 
We prove the latter by induction on $\vt$.

The nontrivial step is the step $\cof\la=\Om$, so that let 
$\vt=\bigcup_{\al<\Om}\vt_\al,$ for an increasing 
$\od$ sequence of ordinals $\vt_\al.$ 
Let $I_\al=[\vt_\al,\vt_{\al+1}).$ 
Then, by the induction hypothesis, for any $\al<\Om$ the set 
$X_\al=\ans{S\res I_\al:S\in X}\sq 2^{I_\al}$ is 
\dd\msl order-isomorphic to an antichain\/ 
$A_\al\sq 2\Lom$ via an\/ $\od$ isomorphism\/ $i_\al,$ 
and the map, which sends\/ $\al$ to\/ $A_\al$ and\/ $i_\al,$ 
is\/ $\od$. 
It follows that the map, which sends each $S\in X$ to the 
concatenation of all sequences $i_\al(x\res I_\al)$, 
is an $\od$ \dd\msl order-isomorphism $X$ onto an antichain 
in $2^\Om.$ 
Therefore it suffices to prove the lemma in the 
case $\vt=\Om.$ 
Thus let $X\sq 2^\Om.$ 

First of all, note that each sequence $S\in X$ is $\rod$.  
Lemma 7 in \cite{ksol}  shows that, in this case, we have 
$S\in\rL[S\res\eta]$ for an ordinal $\eta<\Om.$ 
Let $\eta(S)$ be the least such an ordinal, and 
$h(S)=S\res{\eta(S)},$ so that $h(S)$ is a countable initial 
segment of $S$ and $S\in\rL[h(S)].$ 
Note that $h$ is still $\od$. 

Consider the set $U=\ran h=\ens{h(S)}{S\in X}\sq 2\Lom.$ 
We can assume that every sequence $u\in U$ has a limit length. 
Then $U=\bigcup_{\ga<\Om}U_\ga,$ where 
$U_\ga=U\cap 2^{\om\ga}$ 
($\om\ga$ is the the \dd\ga th limit ordinal). 
For $u\in U_\ga,$ let $\ga_u=\ga$. 

If $u\in U$ then by construction the set $X_u=\ans{S\in X:h(S)=u}$ 
is $\od(u)$ and satisfies $X_u\sq\rL[u]$. 
Therefore, it follows from the known properties of 
the Solovay model that $X_u$ belongs to $\rL[u]$ and 
is of cardinality\/ $\le\Om$ in\/ $\rL[u]$.
Fix an enumeration $X_u=\ans{S_u(\al):\ga_u\le\al<\Om}$ 
for all $u\in U$. 
We can assume that the map $\al,u\longmapsto S_u(\al)$ 
is $\od$. 

If $u\in U$ and $\ga_u\le\al<\Om$, then we define  
a shorter sequence, $s_u(\al)\in 3^{\om\al+1}$, as follows.
\ben
\renu 
\itla{w1}\msur
$s_u(\al)(\xi+1)=S_u(\al)(\xi)$ for any $\xi<\om\al$. 

\itla{w2}\msur
$s_u(\al)(\om\al)=1$.

\itla{w4}
Let $\da<\al.$ 
If $S_u(\al)\res\om\da=S_{v}(\da)\res\om\da$ for some $v\in U$ 
(equal to or different from $u$) then 
$s_u(\al)(\om\da)=0$ whenever $S_u(\al)\lex S_{v}(\da),$ and 
$s_u(\al)(\om\da)=2$ whenever $S_{v}(\da)\msl S_u(\al).$ 

\itla{w3}
Otherwise (\ie, if there is no such $v$),  
$s_u(\al)(\om\da)=1$.
\een
To demonstrate that \ref{w4} is consistent, we show that 
$S_{u'}(\da)\res\om\da=S_{u''}(\da)\res\om\da$ implies $u'=u''.$ 
Indeed, as by definition $u'\subset S_{u'}(\da)$ and  
$u''\subset S_{u''}(\da),$ $u'$ and $u''$ must be 
\dd\sq compatible: let, say, $u'\sq u''.$ 
Now, by definition, $S_{u''}(\da)\in\rL[u''],$ therefore 
$\in \rL[S_{u'}(\da)]$ because 
$u''\sq S_{u''}(\da)\res\om\da=S_{u'}(\da)\res\om\da,$ 
finally $\in \rL[u'],$ which shows that  
$u'=u''$ as $S_{u''}(\da)\in X_{u''}$.

We are going to prove that the map 
$S_u(\al)\longmapsto s_u(\al)$ 
is a \dd\lex order isomorphism, so that 
$S_v(\ba)\lex S_u(\al)$ implies $s_v(\ba)\lex s_u(\al)$. 

We first observe that $s_v(\ba)$ and $s_u(\al)$ are 
\dd\sq incomparable. 
Indeed assume that $\ba<\al.$ 
If $S_u(\al)\res\om\ba\ne S_v(\ba)\res\om\ba$ then clearly 
$s_v(\ba)\not\sq s_u(\al)$ by \ref{w1}. 
If $S_u(\al)\res\om\ba= S_v(\ba)\res\om\ba$ then 
$s_u(\al)(\om\ba)=0\text{ or }2$ by \ref{w4} 
while $s_v(\ba)(\om\ba)=1$ by \ref{w2}. 
Thus all $s_u(\al)$ are mutually \dd\sq incomparable, so that 
it suffices to show that conversely 
$s_v(\ba)\msl s_u(\al)$ implies $S_v(\ba)\msl S_u(\al)$. 
Let $\za$ be the least ordinal such that 
$s_v(\ba)(\za)< s_u(\al)(\za);$ then 
$s_u(\al)\res\za= s_v(\ba)\res\za$ 
and $\za\le\tmin\ans{\om\al,\om\ba}.$ 

The case when $\za=\xi+1$ is clear: then by definition 
$S_u(\al)\res\xi= S_v(\ba)\res\xi$ while 
$S_v(\ba)(\xi)< S_u(\al)(\xi),$ 
so let us suppose that $\za=\om\da,$ where 
$\da\le\tmin\ans{\al,\ba}.$ 
Then obviously $S_u(\al)\res\om\da= S_v(\ba)\res\om\da.$ 
Assume that one of the ordinals $\al,\,\ba$ is equal to 
$\da,$ say, $\ba=\da.$ 
Then $s_v(\ba)(\om\da)=1$ while $s_u(\al)(\om\da)$ is 
computed by \ref{w4}. 
Now, as $s_v(\ba)(\om\da)< s_u(\al)(\om\da)$, we conclude that 
$s_u(\al)(\om\da)=2,$ hence $S_v(\ba)\msl S_u(\al),$ as required. 
Assume now that $\da<\tmin\ans{\al,\ba}.$ 
Then easily $\al$ and $\ba$ appear in one and the same class 
\ref{w4} or \ref{w3} with respect to the $\da$. 
However this cannot be \ref{w3} because 
$s_v(\ba)(\om\da)\ne s_u(\al)(\om\da).$ 
Hence we are in \ref{w4}, so that, for some (unique) 
$w\in U$. 
$0=S_v(\ba)\msl S_w(\da)\lex S_u(\al)=2,$ as required. 

This ends the proof of the lemma, except for the fact 
that the sequences $s_u(\al)$ belong to $3\Lom,$ but improvement 
to $2\Lom$ is easy.
\epf


\punk{Decomposing thin $\od$ sets in the Solovay model}
\las{Td:sm}

Here we prove 
Theorem~\ref{t:sm}.
{\ubf We assume to the contrary} that the $\od$ set 
$\Ua$ of all reals $x\in\Xa$ such that 
there is no $\od$ \dd\cle chain $C$ containing $x$, 
is non-empty.  
If $R\sq\bn$ is an $\od$ set then let $\cfx R$ consist of 
all $\od$ 
maps $F:{\stk{\bn}\cle}\to\stk{A}\lexe$, 
where $A\sq2\Lom$ is an $\od$ antichain, such that  
\ben
\Renu
\itla{F2i}
$F$ is \lr\ order preserving, \ie,
\imar{F2i}
$x\cle y \limp  F(x)\lexe F(y)$ --- in particular, 
$x\apr y \limp  F(x)=F(y)$ --- for all $x,y\in\bn$;

\itla{F2ii}
if $x,y\in R$ are \dd\cle incomparable then $F(x)=F(y)$, 
or equivalently provided \ref{F2i} holds, 
$F(x)\lex F(y)\imp x\cl y$ for all $x,y\in R$. 
\een
We let
$$
x\rE_R y \;\;\text{ iff }\;\; \kaz F\in\cfx R\:(F(x)=F(y))\,. 
$$
Note that a function $F\in\cfx R$ has to be not 
just \dd\apr invariant by \ref{F2i}, 
but also invariant \vrt\ the common equivalence hull of the 
relation $\apr$ 
and the (non-equivalence) relation of 
being \dd\cle incomparable, by \ref{F2ii}. 

Still any ${\rE_R}$ is an $\od$ \eqr.

If $R\sq R'$ then $\cfx{R'}\sq\cfx R$, and hence 
$x\rE_R y$ implies $x\rE_{R'} y$.

\ble
\lam{Lef2}
If\/ $R\sq\bn$ is\/ $\od$, ${\rE_R}\sq H\sq\bn\ti\bn$,
and\/ $H$ is\/ $\od$, 
then there is a function\/ $F\in\cF_R$ such that\/ 
$\kaz x,y\:({F(x)=F(y)}\imp H(x,y))$. 
\ele
\bpf 
Clearly $\card{\cF_R}=\Theta$ and $\cF_R$ admits 
an $\od$ enumeration $\cF_R=\ens{F_\xi}{\xi<\vT}$.
If $x\in\bn$ then let 
$ 
f(x)=F_0(x)\we F_1(x)\we \dots\we F_\xi(x)\we\;\dots
$  
--- the 
concatenation of all sequences $F_\xi(x)$. 
Then $f:\stk{\bn}\cle\to\stk{W}\lexe$ is a LR order 
preserving \od\ map, where 
$W=\ran f=\ens{f(r)}{r\in\bn}\sq2^\vT$, 
and ${f(x)=f(y)}\imp H(x,y)$ by the construction. 
By Lemma~\ref{apal31} there is an \od\ isomorphism 
$g:\stk{W}\lexe\onto\stk{A}\lexe$ onto an antichain 
$A\sq2\Lom.$   
The superposition $F(x)=g(f(x))$ proves the lemma.
\epf

\ble
\lam{SMsmu2}
If\/ $R\sq \Ua$ is a non-empty\/ $\od$ set then the 
condition\/ $R\ti R$ \efa Rs that\/ $\doxl\rE_R\doxr$.\qed
\ele
\bpf
Otherwise, by Lemma \ref{frl}, 
there is a function $F\in\cF_R$ and a condition 
$X\ti Y$ in $\spa R$ with $X\cup Y\sq R$, which \efa Rs 
$F(\doxl)(\xi)=0\ne 1= F(\doxr)(\xi)$ for an   
ordinal $\xi<\Om$. 
We may assume that $X\ti Y$ is a saturated condition. 
Then $F(x)(\xi)=0 \ne 1= F(y)(\xi)$ for any pair 
$\ang{x,y}\in X\ti Y$, so that we have $F(x)\ne F(y)$ and 
$\neg\;(x\er R y)$ 
whenever $\ang{x,y}\in X\ti Y$, which contradicts the choice 
of $X\ti Y$ in $\spa R$.
\epf

\ble
\lam{51sm}
Let\/ $R\sq \Ua$ be a non-empty\/ $\od$ set. 
Then\/ $R\ti R$ 
does {\bfit not}\/ \efa R\ that\/ $\doxl,\doxr$ 
are\/ \dd\cle comparable.
\ele
\bpf
Suppose to the contrary that $R\ti R$ forces the comparability. 
Then by Lemma~\ref{frl} a subcondition $X\ti Y$  
either \efa Rs $\doxl\apr\doxr$ or \efa Rs ${\doxl\cl\doxr}$; 
$X,Y\sq R$ are $\od$ sets and 
$\ek{X}{\rE_R}\cap\ek{Y}{\rE_R}\ne\pu$.\vom

{\ubf Case A:} \ 
condition $X\ti Y$ \efa Rs $\doxl\apr\doxr$. 
We claim that the $\od$ set 
$W=\ens{\ang{x,x'}\in X\ti X}{x\rE_R x'\land x'\napr x}$ 
is empty. 
Indeed otherwise $W$ is a condition in the forcing
$\odd=\ens{P\sq\bn\ti\bn}{\pu\ne P\in\od}$, which is just
the 2-dimensional version of $\OD$ with the same basic
properties.
Note that $\odd$ adds pairs $\ang{\doxl,\doxr}\in W$ of \pge\ 
(separately) reals $\dox,{\dox}{}'\in X$, and the 
condition $W\,$ \dd{\odd}forces that 
${\dox}{}'\rE_R \dox$ and ${\dox}{}'\napr \dox$. 

If $P\in\odw$ then obviously 
$\ek{\dom P}{\rE_R}=\ek{\ran P}{\rE_R}$. 

Consider a more complex forcing 
$\cP$
of all pairs $P\ti Y'$, 
where $P\in\odd$, $P\sq W$, 
$Y'\in\OD$, $Y'\sq Y$, and 
$\ek{\dom P}{\rE_R}\cap\ek{Y'}{\rE_R}\ne\pu$. 
For instance, $W\ti Y\in\cP$.  
Then $\cP$ adds a pair $\ang{\doxl,\doxr}\in W$ and a separate 
real $\dox\in Y$ such that the pairs $\ang{\doxl,\dox}$, 
$\ang{\doxr,\dox}$ belong to $X\ti Y$ and are 
\dd{(\spa R)}generic, 
hence  
$\doxl\apr\dox\apr\doxr$ by the choice of $X\ti Y$. 
On the other hand, $\doxl\napr\doxr$ since the pair belongs 
to $W$, which is a contradiction.

Thus $W=\pu$. 
Then $X$ is a \dd\cle chain: 
indeed if $x,y\in X$ are \dd\cle incomparable 
then by definition we have $x\rE_R y$, hence $x\apr y$, 
contradiction. 
Thus $X$ is an $\od$ \dd\cle chain with $\pu\ne X\sq\Ua$, 
contrary to the definition of $\Ua$. 
\vom

{\ubf Case B:} \ 
condition $X\ti Y$ \efa Rs $\doxl\cl\doxr$. 
We claim that the $\od$ set 
$W'=\ens{\ang{x,y}\in X\ti Y}{x\rE_R y\land x\not\cl y}$ 
is empty. 
Suppose towards the contrary that $W'\ne\pu$. 
Let $X'=\dom W'$.
As $X'\sq R,$ the condition $X'\ti X'$ \efa Rs that 
$\doxl,\doxr$ are \dd\cle comparable.
Therefore there is a condition $A\ti B$ in $\spa R$, 
with $A\cup B\sq X'$, which \efa Rs $\doxl\cl\doxr$; 
for if it forces $\doxr\cl\doxl$ then just consider $B\ti A$ 
instead of $A\ti B$, and it cannot force $\doxr\apr\doxl$ 
by the result in Case A.
Let $Z=\ens{\ang{x,y}\in W'}{x\in A}$.

Consider the forcing notion $\cP$
of all non-empty $\od$ sets
of the form $P\ti B'$, where 
$P\sq Z$, $B'\sq B$, and 
$\ek{B'}{\rE_R}\cap\ek{\dom P}{\rE_R}\ne\pu$ 
(equivalently, $\ek{B'}{\rE_R}\cap\ek{\ran P}{\rE_R}\ne\pu$).
For instance, $Z\ti B\in\cP$. 
Note that $\cP$ adds a pair $\ang{\doxl,\doxr}\in Z$ and
a separate 
real $\dox\in B$ such that both pairs $\ang{\doxl,\dox}$ and 
$\ang{\doxr,\dox}$ are \dd{(\spa R)}generic. 
It follows that $\cP$ forces 
both $\doxl\cl\dox$ (as this pair belongs to $A\ti B$) 
and $\dox\cl\doxr$ (it belongs to $X\ti Y$), 
hence, forces  $\doxl\cl\doxr$. 
On the other hand $\cP$ forces $\doxl\not\cl\doxr$ 
(as this pair belongs to $Z\sq W'$), a contradiction.\vom

Thus $W'=\pu$; 
in other words, if $x\in X,$ $y\in Y$, and $x\rE_R y$ 
then\/ $x\cl y$ strictly.
The $\od$ set 
$\wX=\ens{x'}{\sus x\in X\,(x\rE_R x'\land x'\cle x)}$ 
is downwards\/ \dd\cle closed 
in each\/ \dd{\rE_R}class, $X\sq\wX$, and  
still $Y\cap\wX=\pu$. 

\bcl
If\/ $x\in \wX\cap R$, $y\in R\bez\wX$, and\/ $y\rE_R x$, 
then\/ $x\cl y$. 
\ecl
\bpf
Otherwise, the following $\od$ set 
$$
H_0=\ens{y\in R\bez\wX}
{\sus x\in\wX\cap R\:(x\rE_R y\land x\not\cl y)}
\sq R
$$          
is non-$\pu$. 
As above (Subcase B1), 
there is a condition $H\ti H'$ in $\spa R$, 
with $H\cup H'\sq H_0$, 
which \efa Rs  $\doxl\cl \doxr$, and then,  
by the result in Case B1, $y\cl y'$ holds whenever 
$\ang{y,y'}\in H\ti H'$  and $y\rE_R y'$. 
By construction the $\od$ set 
$$
\wX_1 =\ens{x\in\wX\cap R}
{\sus y'\in H'\,({x\rE_R y'} \land x\not\cl y')} 
$$
satisfies $\ek{\wX_1}{\rE_R}=\ek{H}{\rE_R}=\ek{H'}{\rE_R}$, 
hence $\wX_1\ti H$ is a condition in $\spa R$.
Let $\ang{x_1,y}\in \wX_1\ti H$ be any \dd{(\spa R)}generic pair. 
Then ${x_1\rE_R y}$ by Lemma~\ref{SMsmu2}, 
and, by the choice of $R$ and the result in Case A, 
we have $x_1\cl y$ or $y\cl x_1$. 
Yet by construction $x_1\in\wX$, $y\nin \wX$, and $\wX$ 
is downwards closed in each \dd{\rE_R}class. 
Thus in fact $x_1\cl y$.
Therefore, for all $y'\in H'$, if $x_1\rE_R y'$ then 
$x_1\cl z\cl y'$, which contradicts to $x_1\in\wX_1$.
\epF{Claim} 

We conclude by Lemma~\ref{Lef2} that   
{\it there is a single function\/ $F\in\cF_R$ 
such that if\/ $x\in \wX\cap R$, $y\in R\bez\wX$, and\/ 
$F(x)=F(y)$, then\/ $x\cl y$}.

Prove that {\it the derived  function}
$$
G(x)=\left\{
\bay{rcl}
F(x)\we 0\,,&\text{whenewer}& x\in\wX\\[1ex]
F(x)\we 1\,,&\text{whenewer}& x\in\bn\bez\wX 
\eay
\right.
$$
{\it belongs to\/ $\cfx R.$} 
First of all, still $G\in\cF$ since $C$ is downwards 
\dd\cle closed in each \dd{\rE_R}class. 
Now suppose that $x,y\in R$ and $G(x)\lex G(y)$.
Then either $F(x)\lex F(y)$, or $F(x)=F(y)$ and $x\in\wX$ but 
$y\nin\wX$.
In the ``either'' case immediately $x\cl y$ since $F\in\cfx R$. 
In the ``or'' case we have $x\cl y$ by the choice of $F$ 
and the definition of $G$.
Thus $G\in\cfx R$. 

Now pick any pair of reals $x\in X$ and $y\in Y$ with $x\rE_R y$. 
Then $G(x)=G(y)$ since $G\in\cfx R$. 
But we have $x\in \wX$ and $y\nin\wX$ since 
$X\sq\wX$ and $Y\cap\wX=\pu$ by construction, and in this case 
surely $G(y)\ne G(x)$ by the definition of $G$.
This contradiction completes the proof of Lemma~\ref{51sm}.
\epf

Lemma~\ref{51sm} plus Lemma~\ref{29sm} imply 
Theorem~\ref{t:sm}.\vtm

\qeDD{Theorem~\ref{t:sm}}

\end{document}